\documentclass[10pt]{amsart}

\usepackage{amsmath,amsthm,amsfonts,amscd,amssymb}
\usepackage{hyperref,wasysym}

\newtheorem{thm}{Theorem}
\newtheorem{cor}[thm]{Corollary}
\newtheorem{lem}[thm]{Lemma}

\theoremstyle{definition}

\theoremstyle{remark}

\usepackage[colorinlistoftodos]{todonotes}

\begin{document}

\title{Counting basis extensions in a lattice}

\author{Maxwell Forst}
\author{Lenny Fukshansky}\thanks{Fukshansky was partially supported by the Simons Foundation grant \#519058}

\address{Institute of Mathematical Sciences, Claremont Graduate University, Claremont, CA 91711}
\email{maxwell.forst@cgu.edu}
\address{Department of Mathematics, 850 Columbia Avenue, Claremont McKenna College, Claremont, CA 91711}
\email{lenny@cmc.edu}

\subjclass[2010]{Primary: 11H06, 11C08, 11C20, 11D85, 11D45}
\keywords{integer lattice, basis, extending to a basis, primitive vectors, multilinear form, counting lattice points, Farey fractions}

\begin{abstract}
Given a primitive collection of vectors in the integer lattice, we count the number of ways it can be extended to a basis by vectors with sup-norm bounded by $T$, producing an asymptotic estimate as $T \to \infty$. This problem can be interpreted in terms of unimodular matrices, as well as a representation problem for a class of multilinear forms. In the $2$-dimensional case, this problem is also connected to the distribution of Farey fractions. As an auxiliary lemma we prove a counting estimate for the number of integer lattice points of bounded sup-norm in a hyperplane in~$\mathbb R^n$. Our main result on counting basis extensions also generalizes to arbitrary lattices in~$\mathbb R^n$. Finally, we establish some basic properties of sparse representations of integers by multilinear forms.
\end{abstract}

\maketitle

\def\A{{\mathcal A}}
\def\B{{\mathcal B}}
\def\C{{\mathcal C}}
\def\D{{\mathcal D}}
\def\F{{\mathcal F}}
\def\x{{\mathcal H}}
\def\I{{\mathcal I}}
\def\J{{\mathcal J}}
\def\K{{\mathcal K}}
\def\L{{\mathcal L}}
\def\M{{\mathcal M}}
\def\O{{\mathcal O}}
\def\R{{\mathcal R}}
\def\s{{\mathcal S}}
\def\V{{\mathcal V}}
\def\W{{\mathcal W}}
\def\X{{\mathcal X}}
\def\Y{{\mathcal Y}}
\def\H{{\mathcal H}}
\def\OO{{\mathcal O}}
\def\BB{{\mathbb B}}
\def\cee{{\mathbb C}}
\def\pee{{\mathbb P}}
\def\que{{\mathbb Q}}
\def\real{{\mathbb R}}
\def\zed{{\mathbb Z}}
\def\hyp{{\mathbb H}}
\def\aa{{\mathfrak a}}
\def\qbar{{\overline{\mathbb Q}}}
\def\eps{{\varepsilon}}
\def\ahat{{\hat \alpha}}
\def\bhat{{\hat \beta}}
\def\gt{{\tilde \gamma}}
\def\h{{\tfrac12}}
\def\be{{\boldsymbol e}}
\def\bei{{\boldsymbol e_i}}
\def\bff{{\boldsymbol f}}
\def\ba{{\boldsymbol a}}
\def\bb{{\boldsymbol b}}
\def\bc{{\boldsymbol c}}
\def\bm{{\boldsymbol m}}
\def\bk{{\boldsymbol k}}
\def\bi{{\boldsymbol i}}
\def\bl{{\boldsymbol l}}
\def\bq{{\boldsymbol q}}
\def\bu{{\boldsymbol u}}
\def\bt{{\boldsymbol t}}
\def\bs{{\boldsymbol s}}
\def\bv{{\boldsymbol v}}
\def\bw{{\boldsymbol w}}
\def\bx{{\boldsymbol x}}
\def\bX{{\boldsymbol X}}
\def\bz{{\boldsymbol z}}
\def\bwy{{\boldsymbol y}}
\def\bY{{\boldsymbol Y}}
\def\bL{{\boldsymbol L}}
\def\baa{{\boldsymbol\alpha}}
\def\bbb{{\boldsymbol\beta}}
\def\bet{{\boldsymbol\eta}}
\def\bxi{{\boldsymbol\xi}}
\def\bo{{\boldsymbol 0}}
\def\bol{{\boldkey 1}_L}
\def\ep{\varepsilon}
\def\p{\boldsymbol\varphi}
\def\q{\boldsymbol\psi}
\def\rank{\operatorname{rank}}
\def\aut{\operatorname{Aut}}
\def\lcm{\operatorname{lcm}}
\def\sgn{\operatorname{sgn}}
\def\spn{\operatorname{span}}
\def\md{\operatorname{mod}}
\def\Norm{\operatorname{Norm}}
\def\dim{\operatorname{dim}}
\def\det{\operatorname{det}}
\def\Vol{\operatorname{Vol}}
\def\rk{\operatorname{rk}}
\def\Gal{\operatorname{Gal}}
\def\WR{\operatorname{WR}}
\def\WO{\operatorname{WO}}
\def\GL{\operatorname{GL}}
\def\pr{\operatorname{pr}}
\def\phi{\operatorname{\varphi}}

\section{Introduction}
\label{intro}

A collection of $1 \leq m < n$ vectors $\ba_1,\dots,\ba_m \in \zed^n$ is said to be {\it primitive} if it is extendable to a basis of $\zed^n$, i.e. if there exist some vectors $\bb_1,\dots,\bb_{n-m} \in \zed^n$ such that
$$\ba_1,\dots,\ba_m,\bb_1,\dots,\bb_{n-m}$$
is a basis for $\zed^n$. There is an easy-to-use criterion for primitivity in terms of the $n \times m$ integer matrix 
$$A = (\ba_1\ \dots\ \ba_m),$$
which follows directly from Lemma~2, p.~15, \cite{cassels}: the collection $\ba_1,\dots,\ba_m$ is primitive if and only if its full-rank minors, called {\it Grassmann coordinates} (or Pl\"ucker coordinates), are coprime. This is equivalent to saying that the matrix $A$ is {\it unimodular}, i.e. there exists an $n \times (n-m)$ matrix $B$ such that the matrix $\begin{pmatrix} A & B \end{pmatrix} \in \GL_n(\zed)$. This criterion is also recorded in~\cite{zurich} (Lemma 2), where it is referred to as the PID version of the Quillen-Suslin theorem.

Let $T > 0$ be a real number and write $|\ |$ for the sup-norm on vectors and matrices, i.e. maximum of absolute values of coordinates. Given a primitive collection $\ba_1,\dots,\ba_m \in \zed^n$, we want to understand in how many ways can it be extended to a basis by vectors $\bb_1,\dots,\bb_{n-m}$ with $|\bb_i| \leq T$. Equivalently, given a unimodular $n \times m$ matrix $A$, we want to count the number of integer $n \times (n-m)$ matrices $B$ with $|B| \leq T$ such that $\begin{pmatrix} A & B \end{pmatrix} \in \GL_n(\zed)$. Specifically, we are interested in the asymptotic behavior of this number as $T \to \infty$. We prove the following result.

\begin{thm} \label{main} Let $\ba_1,\dots,\ba_m \in \zed^n$ be a primitive collection of vectors.

\begin{enumerate}

\item If $m < n-1$, the number of vectors $\bb \in \zed^n$ with $|\bb| \leq T$ such that the collection $\ba_1,\dots,\ba_m, \bb$ is again primitive is equal to $\Theta(T^n)$ as $T \to \infty$.

\item If $m = n-1$, the number of vectors $\bb \in \zed^n$ with $|\bb| \leq T$ such that the collection $\ba_1,\dots,\ba_m, \bb$ is a basis for $\zed^n$ is equal to $\Theta(T^{n-1})$ as $T \to \infty$.

\end{enumerate}

\noindent
As a result, for any $1 \leq k < n-m$ there exist $\Theta(T^{nk})$ collections of vectors $\bb_1,\dots,\bb_k \in \zed^n$ with $|\bb_i| \leq T$, $1 \leq i \leq k$, such that $\{ \ba_i, \bb_j : 1 \leq i \leq m, 1 \leq j \leq k \}$ is again primitive. Further, there are $\Theta(T^{n^2-nm-1})$ such collections $\bb_1,\dots,\bb_{n-m}$ so that 
$$\zed^n = \spn_{\zed} \left\{ \ba_1,\dots,\ba_m,\bb_1,\dots,\bb_{n-m} \right\}.$$
The constants in the $\Theta$-notation depend on the vectors $\ba_1,\dots,\ba_m$, $n$ and $m$.
\end{thm}

\noindent
In fact, in part (2) of Theorem~\ref{main} we obtain a more precise asymptotic result with an explicit constant (Corollary~\ref{extend}). We also discuss the constants in $\Theta$-notation for part (1) in Lemma~\ref{add_one}.

More generally, notice that any lattice $\Lambda \subset \real^n$ is of the form $\Lambda = U \zed^n$ for some matrix $U \in \GL_n(\real)$. As such, bases in $\Lambda$ are in bijective correspondence with bases in $\zed^n$, given by multiplication by~$U$. This correspondence allows to extend Theorem~\ref{main} to arbitrary lattices, where we call a collection of vectors in $\Lambda$ primitive if it is a basis or can be extended to a basis of~$\Lambda$.

\begin{cor} \label{gen_lat} Let $\ba_1,\dots,\ba_m$ be a primitive collection of vectors in a full-rank lattice $\Lambda  \subset \real^n$ with $1 \leq m < n$. Then there exist 
$$\Theta \left( T^{n + \min \{ 0, n-m-2 \}} \right)$$
vectors $\bb \in \Lambda$ with $|\bb| \leq T$ so that the collection $\ba_1,\dots,\ba_m,\bb$ is still primitive in~$\Lambda$, and hence for any $1 \leq k < n-m$ there exist $\Theta(T^{nk})$ collections of vectors $\bb_1,\dots,\bb_k \in \Lambda$ with $|\bb_i| \leq T$, $1 \leq i \leq k$, such that $\{ \ba_i, \bb_j : 1 \leq i \leq m, 1 \leq j \leq k \}$ is again primitive. Further, there are $\Theta(T^{n^2-nm-1})$ collections of vectors $\bb_1,\dots,\bb_{n-m} \in \Lambda$ such that $|\bb_i| \leq T$ for each $1 \leq i \leq n-m$ and
$$\Lambda = \spn_{\zed} \left\{ \ba_1,\dots,\ba_m,\bb_1,\dots,\bb_{n-m} \right\}.$$
The constants in the $\Theta$-notation depend on the lattice $\Lambda$, the vectors $\ba_1,\dots,\ba_m$, $n$ and $m$.
\end{cor}

Basis extensions play a key role in reduction theory, which aims to specify a canonical choice of a basis for a lattice, often with some favorable properties, such as short norm and near-orthogonality (see~\cite{CS}, \cite{lek}, \cite{martinet} for details). Some of the best known reduction procedures are Lagrange-Gauss reduction (in two dimensions), HKZ (Hermite-Korkine-Zolotarev) reduction, Minkowski reduction, and LLL (Lenstra-Lenstra-Lov\'asz) reduction, each specifying an algorithm for building a basis with some special properties a vector at a time. For instance, Minkowski reduction procedure starts with a shortest (with respect to a chosen norm) primitive vector $\ba_1 \in \Lambda$, and then on each $i$-th step, $1 < i \leq n$, selects a vector $\ba_i$ which is shortest out of all vectors $\bb$ such that the collection $\ba_1,\dots,\ba_{i-1},\bb$ is primitive. This motivates our counting problem for the number of all possible basis extensions of bounded sup-norm.

Additionally, the distribution properties of unimodular matrices, as well as matrices that can be completed to unimodular have received quite a bit of attention in the recent years (see, for instance, \cite{fang}, \cite{guo}, \cite{zurich}, \cite{zhan}). In particular, in~\cite{zurich} it has been proved that unimodular matrices have natural density among all integer matrices of the same dimensions. In a way, our main result complements this theorem by counting the number of ways a given unimodular matrix can be completed to a matrix in $\GL_n(\zed)$ among all possible choices with sup-norm $\leq T$ as $T \to \infty$.

Another interpretation of our problem is in terms of representation numbers for a particular class of multilinear forms. Indeed, given a primitive collection of vectors $\ba_1,\dots,\ba_m \in \zed^n$ let $A = (\ba_1\ \dots\ \ba_m)$ be the corresponding $n \times m$ unimodular matrix. Let $X = (x_{ij})_{1 \leq i \leq n, 1 \leq j \leq n-m}$ be the $n \times (n-m)$ variable matrix. Then the $n \times n$ matrix $(A\ B)$ is in $\GL_n(\zed)$ if and only if the corresponding multilinear form $\det (A\ X)$ of degree $n-m$ in $n(n-m)$ variables with coprime integer coefficients represents $\pm 1$ at $X = B$. Existence of integer representations by multilinear forms has been recently investigated in \cite{ab-lf}, and in this way the current paper continues that investigation by counting the number of such representations in the special case of determinantal forms.

The paper is organized as follows. In Section~\ref{primitive} we prove part (1) of Theorem~\ref{main} by looking at a direct sum decomposition of $\zed^n$ with respect to the lattice spanned by the columns of the matrix $A$ and counting primitive points modulo this lattice. We prove a counting lemma for the number of integer lattice points in a hyperplane in Section~\ref{lat_hyp}, using a result of Borosh et al.~\cite{borosh} on small-size solutions to linear equations and some previous lattice-point counting estimates of the second author~\cite{lf_int}, \cite{lf_gh}. We then use this lemma in Section~\ref{gln} to finish the proof of Theorem~\ref{main} along with Corollary~\ref{gen_lat}. We also comment on the constants in $\Theta$-notation, obtaining a more precise result at least for part (2) of Theorem~\ref{main}. In Section~\ref{farey} we discuss a connection of our problem in the two-dimensional situation to Farey fractions. Finally, in Section~\ref{sparse} we include some observations on sparse representation of integers by more general multilinear forms in the spirit of the previous paper~\cite{ab-lf}.

\section{Extending a primitive collection of vectors}
\label{primitive}

In this section we count the number of primitive extensions of a primitive collection of $m <n-1$ vectors in $\zed^n$ by one vector, proving the first part of Theorem~\ref{main}. For integers $1 \leq d \leq n$, let us write $[n] := \{ 1,\dots,n \}$ and define the set of indexing subsets
\begin{equation}
\label{index_set}
\J(n,d) := \left\{ I \subseteq [n] : |I| = d \right\},
\end{equation}
then $|\J(n,d)| = \binom{n}{d}$. For a real number $T > 0$, define also the integer $n$-cube centered at the origin with sidelength $2T$ as
$$C_n(T) := \left\{ \bx \in \zed^n : |\bx| \leq T \right\},$$
then $|C_n(T)| = (2[T]+1)^n$. 

\begin{lem} \label{add_one} Let $1 \leq m < n-1$ and let $\ba_1,\dots,\ba_m \in \zed^n$ be a primitive collection of vectors. For $T > 0$, define
$$f(T) = \left| \left\{ \bx \in C_n(T) : \ba_1,\dots,\ba_m,\bx \text{ is primitive} \right\} \right|,$$
then as $T \to \infty$,
\begin{equation}
\label{ft_1}
f(T) \leq \left( \zeta(n)^{-1} + \eps \right) (2T+1)^n,
\end{equation}
where $\zeta$ is the Riemann zeta-function, $\eps > 0$. Additionally,
\begin{equation}
\label{ft_2}
f(T) \geq \beta(n, m, A) T^n,
\end{equation}
where $\beta(n, m, A)$ is a constant depending only on $n, m$ and the matrix $A$.
\end{lem}

\proof
Since the collection $\ba_1,\dots,\ba_m$ is primitive, the corresponding $n \times m$ matrix $A$ is unimodular. By the primitivity criterion stated in Section~\ref{intro}, we want to count $\bx \in C_n(T)$ such that the extended $n \times (m+1)$ matrix $(A\ \bx)$ is still unimodular. First notice that each $\bx$ must itself be a primitive vector, for otherwise the $\gcd$ of its coordinates will be a multiple of all the Grassmann coordinates of the extended matrix $(A\ \bx)$. Therefore the total number of such vectors $\bx$ is no bigger than the number of primitive vectors in $C_n(T)$. It is a well-known fact that the probability of a vector in $C_n(T)$ being primitive is $\zeta(n)^{-1}$ (this result has apparently first been proved by E. Ces\`aro in 1884, but has been re-discovered several times since; see~\cite{zurich} for the references). More specifically, a result of~\cite{nymann} asserts that
$$\left| \left\{ \bx \in C_n(T) : \bx \text{ is primitive} \right\} \right| = \zeta(n)^{-1} T^n + O(T^{n-1}).$$
Taking any $\eps > 0$ then guarantees~\eqref{ft_1} for all sufficiently large $T$.

Now, let $\Lambda := A\zed^m$ be a lattice of rank $m$. Since $A$ is a unimodular matrix, there exists an $n \times (n-m)$ integer matrix $B$ so that the augmented matrix $(A\ B) \in \GL_n(\zed)$. In fact, such $B$ can be chosen so that $\sqrt{\det(B^{\top} B)}$ is bounded by a function of $A$, call it $\alpha(A)$: this can be done, for instance, by a repetitive application of the search-bound presented in Section~1 of~\cite{ab-lf}, after Theorem~1.3. Then 
$$\Omega := B\zed^{n-m} \cong \zed^n/\Lambda \cong \zed^{n-m}$$
is a lattice, and
\begin{equation}
\label{det}
\det(\Omega) = \sqrt{ \det(B^{\top} B)} \leq \alpha(A),\ \det(\Lambda) = \sqrt{ \det(A^{\top} A)}.
\end{equation}
Then $\zed^n = \Omega \oplus \Lambda$, so for every $\bx \in \zed^n$ there exists a unique pair $\bwy \in \Omega$, $\bz \in \Lambda$ such that $\bx = \bwy + \bz$, so
$$|\bx| \leq |\bwy| + |\bz|.$$
Then $\bx$ is such that $(A\ \bx)$ is unimodular if and only if the corresponding $\bwy \in \Omega$ is primitive, i.e. extendable to a basis of $\Omega$. Let $\gamma \in (0,1]$ and notice that
$$g_{\gamma}(T) := \left| \left\{ \bwy + \bz : \bwy \in C_n(\gamma T) \cap \Omega', \bz \in C_n((1-\gamma)T) \cap \Lambda \right\} \right| \leq f(T),$$
where $\Omega'$ stands for the set of primitive points in $\Omega$. Now notice that 
\begin{equation}
\label{g-eq-1}
g_{\gamma}(T) = \left| C_n(\gamma T) \cap \Omega' \right| \cdot \left| C_n((1-\gamma)T) \cap \Lambda \right|.
\end{equation}
Assume $T \geq \max \left\{ \frac{(n-m) \det(\Omega)}{2 \gamma}, \frac{m \det(\Lambda)}{2(1-\gamma)} \right\}$. Then Lemma 3.1 of~\cite{lf_gh} guarantees that
\begin{equation}
\label{g-eq-2}
\left| C_n(\gamma T) \cap \Omega \right| \geq \left( \frac{2 \gamma T}{(n-m) \det(\Omega)} - 1 \right) \left( \frac{2 \gamma T}{n-m} - 1 \right)^{n-m-1},
\end{equation}
\begin{equation}
\label{g-eq-3}
\left| C_n((1-\gamma) T) \cap \Lambda \right| \geq \left( \frac{2 (1-\gamma) T}{m \det(\Lambda)} - 1 \right) \left( \frac{2 (1-\gamma) T}{m} - 1 \right)^{m-1}.
\end{equation}
Again, by Ces\`aro's theorem the proportion of primitive points among all points in $\Omega$ is $\zeta(n-m)^{-1}$. Combining this observation with \eqref{g-eq-1}, \eqref{g-eq-2}, \eqref{g-eq-3}, \eqref{det} and taking $\gamma = 1/2$, we obtain
$$f(T) \geq g_{1/2}(T) \geq \beta(n, m, A) T^n$$
for an appropriate constant $\beta(n, m, A)$. This proves~\eqref{ft_2}.
\endproof


\section{Counting lattice points in a hyperplane}
\label{lat_hyp}

In this section we prove a counting lemma on the number of integer lattice points in a section of the cube $C_n(T)$ by a hyperplane, building on a previous result for a section by a subspace. Let 
$$L(x_1,\dots,x_n) = \sum_{i=1}^n c_i x_i \in \zed[x_1,\dots,x_n]$$
be a linear form in $n \geq 2$ variables with coprime coefficients, and write $\bc = (c_1,\dots,c_n)$ for this coefficient vector. Let $b \in \zed$ and let $T > 0$ be a real number. Define the set 
$$\BB_{L,b}(T) = \left\{ \bx \in \zed^n : L(\bx) = b,\ |\bx| \leq T \right\} = C_n(T) \cap \left\{ \bx \in \zed^n : L(\bx) = b \right\}.$$
Since coefficients of $L$ are coprime, the equation $L(\bz) = b$ has infinitely many integer solutions for any $b \in \zed$, and so the set $\BB_{L,b}(T)$ is not empty for a sufficiently large $T$. We want to estimate the size of $\BB_{L,b}(T)$ as a function of the coefficients of $L$, $b$ and~$T$. One tool that we will need for this is the Brill-Gordan duality principle (see~\cite{gordan}, as well as Theorem~1 on p.~294 of~\cite{hodge_pedoe}; see also proof of Theorem~4.3 of~\cite{lf_int}, as well as~\cite{heath-brown}, \cite{bertrand} for more contemporary accounts of this principle).

\begin{lem} [Duality Principle] \label{duality} Let $1 \leq m < n$, and let $A, B$ be respectively $n \times m$ and $(n-m) \times n$ integer matrices such that
$$A \zed^m = \left\{ \bx \in \zed^n : B \bx = \bo \right\}.$$
Write $\Delta_{i_1,\dots,i_m}$ for the Grassmann coordinate of $A$, which is the determinant of the submatrix of $A$ whose rows are indexed by $i_1,\dots,i_m \in \{1,\dots,n\}$. Write $\delta^{j_1,\dots,j_{n-m}}$ for the Grassmann coordinate of $B$, which is the determinant of the submatrix of $B$ whose columns are indexed by $j_1,\dots,j_{n-m} \in \{1,\dots,n\}$. Then
$$\Delta_{i_1,\dots,i_m} = (-1)^{i_1+\cdots+i_m}\ \gamma\ \delta^{i_{m+1},\dots,i_n}$$
for an appropriate constant $\gamma \in \que$, where $\{i_1,\dots,i_m,i_{m+1},...,i_n\} = \{1,\dots,n\}$. If column vectors of $A$ and row vectors of $B$ are primitive (i.e., can be extended to a basis of~$\zed^n$), then $\gamma = 1$.
\end{lem}

\begin{thm} \label{lin_form} For any $T$,
\begin{equation}
\label{L_upper}
|\BB_{L,b}(T)| \leq \left( \frac{2(T + \max \{ |L|, |b| \})}{|L|} + 1 \right) (2(T + \max \{ |L|, |b| \} ) + 1)^{n-2},
\end{equation}
and for $T \geq \max \{ |L|, |b| \}$
\begin{equation}
\label{L_lower}
|\BB_{L,b}(T)| \geq \left( \frac{2(T - \max \{ |L|, |b| \})}{|L|} - 1 \right) (2(T - \max \{ |L|, |b| \} ) - 1)^{n-2}.
\end{equation}
Therefore
\begin{equation}
\label{L_O}
|\BB_{L,b}(T)| \sim \frac{(2T)^{n-1}}{|L|}
\end{equation}
as $T \to \infty$.
\end{thm}

\proof
Let
$$\Lambda_L = \left\{ \bx \in \zed^n : L(\bx) = 0 \right\} = \left\{ \bx \in \zed^n : \bc \cdot \bx = 0 \right\},$$
then $\Lambda_L$ is a sublattice of $\zed^n$ of rank $n-1$. Further, if $A$ is any basis matrix for $\Lambda_L$ then column vectors of $A$ must be primitive, since $\Lambda_L$ is the full intersection of $\zed^n$ with a subspace. We define the Grassmann coordinates of $\Lambda_L$ to be the absolute values of Grassmann coordinates of $A$. This definition does not depend on the choice of a basis matrix for $\Lambda_L$, since for any two such basis matrices $A_1,A_2$ there exists a matrix $U \in \GL_{n-1}(\zed)$ such that $A_2 = U A_1$, where $\det(U) = \pm 1$. Let us write $\Delta$ for the maximum of Grassmann coordinates of $\Lambda_L$, then by Lemma~\ref{duality},
\begin{equation}
\label{Delta_L}
\Delta = |\bc| = |L|.
\end{equation}

For a fixed integer $b$, let
$$\Lambda_L(b) = \left\{ \bx \in \zed^n : L(\bx) = b \right\},$$
so $\Lambda_L = \Lambda_L(0)$ and $\BB_{L,b}(T) = \left\{ \bx \in \Lambda_L(b) : |\bx| \leq T \right\}$. Pick any $\bz \in \Lambda_L(b)$, then it is easy to notice that
\begin{equation}
\label{Lb}
\Lambda_L(b) = \left\{ \bx + \bz : \bx \in \Lambda_L \right\},
\end{equation}
i.e. $\bx \mapsto \bx+\bz$ is a bijective map between $\Lambda_L$ and $\Lambda_L(b)$ for any fixed $\bz \in \Lambda_L(b)$. The main theorem of~\cite{borosh} guarantees that there exists $\bz \in \Lambda_L(b)$ such that
\begin{equation}
\label{borosh_bound}
|\bz| \leq \max \{ |L|, |b| \},
\end{equation}
so from now on we use description \eqref{Lb} for $\Lambda_L(b)$ with $\bz$ satisfying \eqref{borosh_bound}. Hence for any $\bwy = \bx + \bz \in \Lambda_L(b)$,
$$|\bwy| \leq |\bx| + |\bz| \leq |\bx| + \max \{ |L|, |b| \}.$$
Combining Theorem 4.2 of~\cite{lf_int} with~\eqref{Delta_L}, we have
\begin{equation}
\label{B_upper}
|\BB_{L,0}(T)| \leq \left( \frac{2T}{|L|} + 1 \right) (2T+1)^{n-2},
\end{equation}
and combining Lemma 3.1 of~\cite{lf_gh} (see also equation (50)) with~\eqref{Delta_L}, we have for every $T \geq \frac{|L|}{2}$
\begin{equation}
\label{B_lower}
|\BB_{L,0}(T)| \geq \left( \frac{2T}{|L|} - 1 \right) (2T-1)^{n-2}.
\end{equation}

Suppose that $\bwy \in \BB_{L,b}(T)$, then $|\bwy| \leq T$ and $\bwy = \bx+\bz$ for a unique $\bx \in \Lambda_L$, so 
$$|\bx| = |\bwy - \bz| \leq T + |\bz| = T + \max \{ |L|, |b| \}.$$
Therefore
$$|\BB_{L,b}(T)| \leq |\BB_{L,0}(T + \max \{ |L|, |b| \})|,$$
and combining this observation with~\eqref{B_upper}, we obtain~\eqref{L_upper}.

Next assume $\bx \in \BB_{L,0}(T -  \max \{ |L|, |b| \})$, which implicitly implies that $T \geq \max \{ |L|, |b| \}$. Let $\bwy = \bx + \bz \in \BB_{L,b}(T)$, then
$$|\bwy| \leq |\bx| + |\bz| \leq T,$$
and so
$$|\BB_{L,b}(T)| \geq |\BB_{L,0}(T - \max \{ |L|, |b| \})|.$$
Then combining this observation with~\eqref{B_lower}, we obtain~\eqref{L_lower}, since we have $T \geq \max \{ |L|, |b| \} > \frac{|L|}{2}$.

Now notice that both, the upper bound~\eqref{L_upper} and the lower bound~\eqref{L_lower} when expanded under the assumption $T \to \infty$ have the order of magnitude $\frac{(2T)^{n-1}}{|L|} + o(T^{n-1})$. Thus
$$\lim_{T \to \infty} \frac{|\BB_{L,b}(T)|}{(2T)^{n-1} / |L|} = 1,$$
which implies~\eqref{L_O}.
\endproof


\section{Extending a collection of vectors to a basis}
\label{gln}

In this section we prove Theorem~\ref{main} and Corollary~\ref{gen_lat}. To start with, let $n \geq 2$, $\ba_1,\dots,\ba_{n-1} \in \zed^n$ a primitive collection of vectors, and let $A = (\ba_1\ \dots \ \ba_{n-1})$ be the corresponding $n \times (n-1)$ unimodular matrix. In how many ways can this primitive collection be extended to a basis of $\zed^n$? More precisely, for a positive integer $T$ let
\begin{equation}
\label{BT}
\BB_A(T) = \left\{ \bz \in \zed^n : \zed^n = \spn_{\zed} \{ \ba_1,\dots,\ba_{n-1}, \bz \}, |\bz| \leq T \right\}.
\end{equation}
We want to understand how big is the cardinality of this set, $|\BB_A(T)|$ as a function of $A$ and $T$. Notice that $\bz \in \BB_A(T)$ if and only if $|\bz| \leq T$ and
$$\det (A\ \bz) = \pm 1.$$
For each $1 \leq k \leq n$, let $A_k$ be the $(n-1) \times (n-1)$ submatrix of $A$ obtained by deleting $k$-th row, then
$$L_A(\bz) := \det (A\ \bz) = \sum_{k=1}^n (-1)^{n+k} \det(A_k) z_k,$$
which is a linear form in the variables $z_1,\dots,z_n$. Since the collection of vectors $\ba_1,\dots,\ba_{n-1}$ is extendable to a basis for~$\zed^n$, it must be true that
$$\gcd \left(\det(A_1),\dots,\det(A_n) \right) = 1,$$
and hence the equation $L_A(\bz) = \pm 1$ has infinitely many integer solutions. Define
$$\Delta_A := \max \{ |\det(A_k)| : 1 \leq k \leq n \},$$
then $|L_A| = \Delta_A \geq 1$, and so we can apply Theorem~\ref{lin_form} with $b=1$ and with $b=-1$ to obtain the following bound.

\begin{cor} \label{extend} For any $T$,
$$|\BB_A(T)| \leq 2 \left( \frac{2T}{\Delta_A} + 3 \right) (2(T + \Delta_A ) + 1)^{n-2},$$
and for $T \geq \Delta_A$
$$|\BB_A(T)| \geq 2 \left( \frac{2T}{\Delta_A} - 3 \right) (2(T - \Delta_A ) - 1)^{n-2}.$$
Therefore
$$|\BB_A(T)| \sim 2 \left( \frac{(2T)^{n-1}}{\Delta_A} \right)$$
as $T \to \infty$.
\end{cor}

\proof
Since
$$\BB_A(T) = \BB_{L_A,1}(T) \cup \BB_{L_A,-1}(T),$$
we are applying Theorem~\ref{lin_form} twice, with $b = \pm 1$, and adding the results. This produces the factor of two in our bounds.
\endproof

Now we combine Corollary~\ref{extend} with Lemma~\ref{add_one} to prove Theorem~\ref{main}.

\proof[Proof of Theorem~\ref{main}]
Parts (1) and (2) of the theorem are given by Lemma~\ref{add_one} and Corollary~\ref{extend}, respectively. Let us prove that there exist $\Theta(T^{n^2-nm-1})$ collections of vectors $\bb_1,\dots,\bb_{n-m} \in \zed^n$ such that $|\bb_i| \leq T$ for each $1 \leq i \leq n-m$ and $\left\{ \ba_1,\dots,\ba_m,\bb_1,\dots,\bb_{n-m} \right\}$ is a basis for~$\zed^n$. 

Let us argue by induction on $n-m \geq 1$. If $n-m=1$, then we only need to add one vector to this primitive collection, and by Corollary~\ref{extend} there are $\Theta(T^{n-1})$ ways to do it. Notice that in this case
$$n^2-nm-1 = n^2 - n(n-1) - 1 = n-1,$$
so the result follows.

Then assume $n-m > 1$, and result is proved for $n-m-1$, i.e. for a primitive collection of $m+1$ vectors. By Lemma~\ref{add_one} there are $\Theta(T^n)$ to extend this primitive collection by one vector. For each such vector, there are $\Theta(T^{n^2-n(m+1)-1})$ extensions to a basis by the induction hypothesis, and hence the total number of extensions of our primitive collection is
$$\Theta(T^n T^{n^2-n(m+1)-1}) = \Theta(T^{n^2-nm-n-1+n}) = \Theta(T^{n^2-nm-1}).$$

Finally, the argument for extending the primitive collection $\{ \ba_1,\dots,\ba_m \}$ to a primitive collection $\{ \ba_i, \bb_j : 1 \leq i \leq m, 1 \leq j \leq k\}$, $1 \leq k < n-m$, is the same as above, but simpler: we do not need to account for the case of the last vector contributing only $\Theta(T^{n-1})$ possibilities, and hence the total number is simply $\Theta(T^{nk})$. This completes the proof.
\endproof

We now extend these observations to general lattices.

\proof[Proof of Corollary~\ref{gen_lat}]
Let $\Lambda$ be a lattice of full rank in $\real^n$, and let 
$$\lambda_1(\Lambda) := \min \left\{ |\bz| : \bz \in \Lambda \setminus \{ \bo \} \right\}$$
be the first successive minimum of $\Lambda$ with respect to the sup-norm. By Minkowski reduction (see, for instance, Theorem~2 on p.66 of \cite{lek} combined with Theorem~2 on p.62 of the same book), there exists a basis $\bz_1,\dots,\bz_n$ for $\Lambda$ such that
$$\frac{1}{n!} \det (\Lambda) \leq \prod_{i=1}^n |\bz_i| \leq \left( \frac{3}{2} \right)^{\frac{(n-1)(n-2)}{2}} \det (\Lambda).$$
Let $U$ be the basis matrix for $\Lambda$ with column vectors $\bz_1,\dots,\bz_n$, ordered in order of increasing sup-norm, so $|U| = |\bz_n|$, and thus
\begin{equation}
\label{U_bnd}
\left( \frac{\det (\Lambda)}{n!} \right)^{1/n} \leq |U| \leq \left( \frac{3}{2} \right)^{\frac{(n-1)(n-2)}{2}} \frac{\det (\Lambda)}{\lambda_1(\Lambda)^{n-1}}.
\end{equation}
Let $\ba_1,\dots,\ba_m$ be a primitive collection of vectors in $\Lambda$, $1 \leq m < n$. Then for each $1 \leq i \leq m$, $\ba_i = U \ba'_i$ for some $\ba'_i \in \zed^n$. Let us write $A = (\ba_1\ \dots\ \ba_m)$, then there exists an $n \times (n-m)$ matrix $B$ such that
$$(A\ B) \zed^n = \Lambda = U \zed^n,$$
hence $U^{-1} (A\ B) = ( (U^{-1}A)\ (U^{-1} B)) \in \GL_n(\zed)$, where $\ba'_1,\dots,\ba'_m$ are the column vectors of $A' := U^{-1} A$. This means that the collection of vectors $\ba'_1,\dots,\ba'_m$ is primitive in $\zed^n$, and hence we can apply Theorem~\ref{main} to it. 

By analogy with~\eqref{BT}, let
$$\BB^m_{A',\zed^n}(T) = \left\{ \bz \in \zed^n : \ba_1,\dots,\ba_m, \bz \text{ is primitive in } \zed^n, |\bz| \leq T \right\},$$
$$\BB^m_{A,\Lambda}(T) = \left\{ \bz \in \Lambda : \ba_1,\dots,\ba_m, \bz \text{ is primitive in } \Lambda, |\bz| \leq T \right\}.$$
Suppose that $\bb \in \BB_{A,\Lambda}(T)$. Then $\ba'_1,\dots,\ba'_m, \bb'$ is primitive in $\zed^n$, where $\bb' = U^{-1} \bb$, and so
$$|\bb'| \leq n |U^{-1}| |\bb| \leq n |U^{-1}| T.$$
Therefore
\begin{equation}
\label{gen1}
| \BB^m_{A,\Lambda}(T) | \leq | \BB^m_{A',\zed^n}(n |U^{-1}| T) | = \Theta \left( T^{n + \min \{ 0, n-m-2 \}} \right),
\end{equation}
since $\min \{ 0, n-m-2 \} = 0$ if $m < n-1$ and $\min \{ 0, n-m-2 \} = -1$ if $m = n-1$.
On the other hand, assume that $\bb' \in \BB_{A',\zed^n}(T/n|U|)$. Then $\ba_1,\dots,\ba_m, \bb'$ is primitive in $\Lambda$, where $\bb = U \bb'$, and so
$$|\bb| \leq n |U| |\bb'| \leq T.$$
Therefore
\begin{equation}
\label{gen2}
| \BB^m_{A,\Lambda}(T) | \geq | \BB^m_{A',\zed^n}(T/n|U|) | = \Theta \left( T^{n + \min \{ 0, n-m-2 \}} \right).
\end{equation}
Combining \eqref{gen1} and \eqref{gen2} and applying an argument identical to the one in the proof of Theorem~\ref{main} above yields the corollary. Since we choose $U$ to be a Minkowski reduced basis for $\Lambda$ with sup-norm bounded as in~\eqref{U_bnd}, the constants in $\Theta$-notation depend intrinsically on $\Lambda$, not on the choice of a basis for~$\Lambda$.
\endproof


\section{Farey fractions and bases in two dimensions}
\label{farey}

In this section we focus on the 2-dimensional case of the problem considered in Section~\ref{gln}: given a primitive vector in $(a,b) \in \zed^2$, in how many ways can it be extended to a basis of $\zed^2$ by a vector $(z_1,z_2)$ of sup-norm $\leq T$? This is equivalent to counting the number of integer solutions to
\begin{equation}
\label{ab_2}
az_2 - bz_1 = \pm 1
\end{equation}
with $|z_1|,|z_2| \leq T$, i.e. the number of points in $\BB_A(T)$ where $A = (a\ b)$. Applying Corollary~\ref{extend}, we have
\begin{equation}
\label{2-dim_ineq}
\frac{4T}{|A|} - 6 \leq |\BB_A(T)| \leq \frac{4T}{|A|} + 6,
\end{equation}
where $|A| = \max \{ |a|, |b| \}$. Here we do not prove any new results, but instead show a connection of this problem to Farey fractions and Diophantine approximation.

The set of rational numbers in the interval $[0,1]$ can be organized into Farey series as follows. For each $n \geq 1$, let $\F_n$ be the set of all rationals $a/b \in [0,1]$ with $\gcd(a,b) = 1$ and $b \leq n$ written in ascending order. For example,
$$\F_5 = \left\{ \frac{0}{1}, \frac{1}{5}, \frac{1}{4}, \frac{1}{3}, \frac{2}{5}, \frac{1}{2}, \frac{3}{5}, \frac{2}{3}, \frac{3}{4}, \frac{4}{5}, \frac{1}{1} \right\}.$$
The set $\F_n$ is called the {\it Farey series of order $n$}. The set $\que \cap [0,1]$ can then be viewed as the limit of $\F_n$ as $n \to \infty$, and this interpretation induces one possible enumeration on $\que \cap [0,1]$. A good source of information on Farey series is Chapter~3 of Hardy and Wright's classical book~\cite{hardy_wright}.

On the other hand, reduced fractions correspond to primitive integer points in the plane. Let
$$\zed_{\pr}^2 = \left\{ (x,y) \in \zed^2 : \gcd(x,y) =1 \right\}.$$
Elements of this set are precisely primitive vectors in $\zed^2$, sometimes also called {\it visible} lattice points, the second name alluding to the property that the line segment connecting $(x,y)$ to the origin contains no other lattice points on it, so $(x,y)$ is not obstructed by anything, hence visible from the origin. If a pair of vectors $\bx_1,\bx_2 \in \zed^2$ forms a basis for the lattice $\zed^2$, then they both must be contained in $\zed_{\pr}^2$ (we routinely identify vectors with their endpoints).

\begin{lem} \label{farey_mink} Let $\bx_1 = (a,b)$ and $\bx_2 = (c,d)$ be in $\zed^2_{\pr}$ and let $n = \max \{ b,d \}$. Then $\bx_1,\bx_2$ form a basis for $\zed^2$ if and only if $\frac{a}{b}$ and $\frac{c}{d}$ are consecutive elements in the Farey series $\F_n$; we call such elements Farey neighbors.
\end{lem}

\proof
First notice that $\bx_1,\bx_2$ form a basis for $\zed^2$ if and only if
$$\left| \det \begin{pmatrix} a & c \\ b& d \end{pmatrix} \right| = \left| ad - bc \right| = 1.$$
Now, suppose that $\frac{a}{b}$ and $\frac{c}{d}$ are Farey neighbors in the Farey series $\F_n$. Then Theorem~28 of~\cite{hardy_wright} guarantees that
$$bc - ad = 1,$$
and so $\bx_1,\bx_2$ are a basis for~$\zed^2$.

In the reverse direction, assume $\bx_1,\bx_2$ are a basis for~$\zed^2$. Assume without loss of generality that $\frac{a}{b} < \frac{c}{d}$. Then $\frac{a}{b}, \frac{c}{d} \in \F_n$, and we only need to prove that there does not exist some $\frac{h}{k} \in \F_n$ such that
\begin{equation}
\label{hk}
\frac{a}{b} < \frac{h}{k} < \frac{c}{d}.
\end{equation}
Let $P$ be the parallelogram spanned by the vectors $\bx_1,\bx_2$, then the vertices of $P$ are $(0,0), (a,b), (c,d), (a+c,b+d)$ and the area of $P$ is the determinant $bc - ad = 1$. Further, $P$ does not contain any integer lattice points in its interior, in particular $(a+c,b+d)$ is also a primitive lattice point. But since $b+d > n$, a primitive point $(h,k)$ satisfying~\eqref{hk} would have to be in the interior of~$P$, hence such a point cannot exist. This proves the lemma.
\endproof

Let
$$C(T) = \left\{ \bz \in \zed_{\pr}^2 : |\bz| \leq T \right\}.$$
We can subdivide $C(T)$ into eight pieces $Q_i^{\pm}(T)$, where $1 \leq i \leq 4$ indicates a quadrant (numbered in the counterclockwise order) and $\pm$ indicates whether the region is above or below the corresponding line $y = \pm x$. For instance,
$$Q^+_1(T) = \left\{ \bz \in C(T) : 0 \leq z_1 \leq z_2 \right\},\ Q_2^-(T) =  \left\{ \bz \in C(T) : 0 \leq z_2 \leq -z_1 \right\}.$$
These pieces have equal cardinality, since they can be obtained from each other by an appropriate reflection. For instance
$$-Q^+_1(T) = \left\{ -\bz \in C(T) : 0 \leq z_1 \leq z_2 \right\} = \left\{ \bz \in C(T) : 0 \leq -z_1 \leq -z_2 \right\} = Q_3^-(T).$$
It is then easy to see that if some $(a,b) \in Q_i^{\pm}(T)$, then all the corresponding vectors extending $(a,b)$ to a basis of $\zed^2$ are contained in $\pm Q_i^{\pm}(T)$: this follows from~\eqref{ab_2}. In other words, $\BB_A(T) \subseteq Q_i^{\pm}(T) \cup - Q_i^{\pm}(T)$, and $|\BB_A(T) \cap Q_i^{\pm}(T)| = |\BB_A(T) \cap -Q_i^{\pm}(T)|$. Further, these cardinalities do not depend on which $Q_i^{\pm}(T)$ the vector $(a,b)$ belongs to. Hence we can assume that $(a,b) \in Q_1^+(T)$, so $|\BB_A(T)| = 2 |\BB_A(T) \cap Q_1^+(T)|$. Then the fraction $a/b$ belongs to the Farey series $\F_n$ for every $n \geq b$. Further, in this case (assuming $T \geq b$)
$$|\BB_A(T) \cap Q_1^+(T)| = \left| \left\{ c/d \in \F_n : b \leq n \leq T, a/b \text{ and } c/d \text{ are neighbors in } \F_n \right\} \right|.$$
Assume that $a/b$ and $c/d$ are neighbors in some $\F_n$, then $n < b+d$ (Theorem~30 of \cite{hardy_wright}) and the next neighbor that will ``squeeze in" between $\frac{a}{b}, \frac{c}{d}$ will be $\frac{a+c}{b+d}$ (Theorem~29 of \cite{hardy_wright}). When $T \gg b$, new neighbors will appear every time $n$ grows by another $b$, and on this interval in $n$, say $(k-1)b \leq n \leq kb$ for some $k$, $a/b$ will acquire two new neighbors: on the left and on the right. This means that 
$$|\BB_A(T) \cap Q_1^+(T)| \sim \frac{2T}{b},$$
and hence $|\BB_A(T)| \sim \frac{4T}{b}$ as $T \to \infty$. Since $a \leq b = |A|$, this agrees with~\eqref{2-dim_ineq}, and also implies that the number of Farey neighbors of a given Farey fraction grows linearly with the denominator.

Farey fractions are also related to Diophantine approximations. Dirichlet's approximation theorem guaranties that for any irrational $\alpha \in \real$ there exist infinitely many primitive points $(p,q) \in \zed^2$ such that 
\begin{equation}
\label{dirichlet}
\left| \alpha - \frac{p}{q} \right| \leq \frac{1}{q^2}.
\end{equation}
Let
$$\D_n(\alpha) = \left\{ p/q \in \que : p/q \text{ satisfies \eqref{dirichlet}}, q \leq n \right\}$$
be the set of all Dirichlet approximations to $\alpha$ with denominator no bigger than~$n$. Farey fractions provide another method of approximating irrational number in the interval $(0,1)$. Let $0 < \alpha < 1$ be irrational, and define the sequence of Farey approximations for $\alpha$ in the following manner: $F_0(\alpha) = \frac{0}{1}$, $F_1(\alpha) = \frac{1}{1}$ and for each $k \geq 2$, $F_k(\alpha) = \frac{a+c}{b+d}$, where
$$\frac{a}{b} = \min_{0 \leq j <k} \left\{ F_j(\alpha): F_j(\alpha) > \alpha \right\},\ \frac{c}{d} = \max_{0 \leq j < k} \left\{ F_j(\alpha): F_j(\alpha) < \alpha \right\},$$
and $\gcd(a,b) = \gcd(c,d) =1$. Define
$$\F_n(\alpha) = \F_n \cap \{ F_k(\alpha) \}_{k=0}^n$$
 to be the set of all Farey approximations to $\alpha$ with denominator~$\leq n$. An element $a/b$ of $\F_n(\alpha)$ is not guarantied to satisfy \eqref{dirichlet}, but is the best upper or lower approximation to $\alpha$ with denominator $\leq b$. Moreover, if $c/d \in \F_n$ is not a Farey approximation, then there exists $a/b \in \F_n$ such that either $c/d < a/b < \alpha$ or $\alpha < a/b < c/d$. Since $b,d \leq n$,
$$\left|\alpha - \frac{a}{b}\right| > \left|\frac{c}{d}- \frac{a}{b}\right|\geq \frac{1}{n^2}.$$
Therefore
\begin{equation}
\label{dir-far}
\D_n(\alpha) \subseteq \F_n(\alpha).
\end{equation}
Now, let $\alpha = [a_0; a_1, a_2, \dots]$ be the continued fraction expansion for $\alpha$, and for each $n \geq 1$ let $\alpha_n = [a_0 ; a_1, a_2, \dots, a_n]$ be its $n$-th convergent. It is well known that
\begin{equation}
\label{cont-dir}
\{ \alpha_k \}_{k=1}^n \subseteq \D_n(\alpha),
\end{equation}
and the convergents alternate in the following sense: $\alpha_{k-1} < \alpha \Leftrightarrow \alpha_k  > \alpha$. We can now prove that, unlike the number of Farey neighbors of a given Farey fraction, the number of Farey approximations of a given irrational number grows less than linearly with the denominator.

\begin{lem} \label{farey_app} Let $0<\alpha<1$, $\alpha \not \in \que$. Then 
$$\lim_{n\rightarrow \infty} \frac{|\F_n(\alpha)| }{n} = 0.$$
\end{lem}

\proof
Let $d_k$ be the denominator of $F_k(\alpha)$ expressed in lowest terms, where $d_1=1, d_2=1$ corresponding to $\frac{0}{1}, \frac{1}{1} \in \F_1(\alpha)$, respectively. Define $a_1=1, b_1=1$, then $d_{k+1} = a_k + b_k$, and 
$$a_{k+1} = d_{k+1},\ b_{k+1} = b_k,$$
if $F_k(\alpha) > \alpha$, or
$$a_{k+1} = a_k,\ b_{k+1} = d_{k+1},$$
if $F_k(\alpha) < \alpha$, where $a_k$ is the denominator of $\min_{0 \leq j < k} \{F_j(\alpha): F_j(\alpha)>\alpha\}$ and $b_k$ is the denominator of $\max_{0\leq j < k} \{F_j(\alpha): F_j(\alpha)<\alpha\}$. Then observe that, with the exception of the first and second terms, the sequence $\{d_k\}$ is strictly increasing and the sequences $\{a_k\}, \{b_k\}$ are non-decreasing. 
Observe also that
$$|\F_n(\alpha)| = |\{d_k:d_k \leq n\}| = \max\{k:d_k \leq n\}.$$
For a fixed $N \in \zed_{>0}$, let $l = |\{d_k: d_k \leq N\}|$ and notice that for $k > l$
$$d_k = a_{k-1}+b_{k-1} \geq d_{k-1} + \min \{a_{k-1},b_{k-1}\}\geq d_l + (k-l)\min \{a_l,b_l\}.$$
This implies 
\begin{eqnarray}
\label{bound}
           \lim_{n\rightarrow\infty} \frac{|\{d_k:d_k \leq n\}|}{n}=& \lim_{n\rightarrow \infty}\left(\frac{|\{d_k: d_k \leq N\}|}{n} + \frac{|\{d_k: N <d_k \leq n\}|}{n}\right) \nonumber \\
 \leq&  \lim_{n\rightarrow \infty} \left( \frac{l}{n}+\frac{\frac{n-N}{\min \{a_l,b_l\}}}{n} \right) = \frac{1}{\min\{a_l,b_l\}}.
\end{eqnarray}
It remains to show that $\{a_k\}, \{b_k\}$ are unbounded. Observe that $a_k$ increases whenever $F_k(\alpha) > \alpha$ and $b_k$ increases whenever $F_k(\alpha) < \alpha$, so we must show that $F_k(\alpha)$ ``switches sides" sufficiently often. By \eqref{cont-dir} and \eqref{dir-far} we know that there exists a sub-sequence $\{k_j\}$ such that $\alpha_j = F_{k_j}(\alpha)$. Now,
$$F_{k_j}(\alpha) >\alpha \Rightarrow F_{k_{j+1}}(\alpha) <\alpha \Rightarrow F_{k_{j+2}}(\alpha) >\alpha,$$
since the continued fraction convergents alternate. By the recurrence relation on $a_k$ and $b_k$, if $F_{k_j}(\alpha) > \alpha$ then
$$b_{k_{j+1}} \geq (k_{j+1}-k_j)a_{k_j}+b_{k_j}\geq f_j,$$
where $f_j$ is the $j$-th Fibonacci number. Likewise, if $F_{k_j}(\alpha) <\alpha$ then
$$a_{k_{j+1}}\geq (k_{j+1}-k_j)b_{k_j}+a_{k_j}\geq f_j.$$
Combining these observations with \eqref{bound}, we conclude that
$$\lim_{n \rightarrow \infty} \frac{\max \{k: d_k \leq n\} }{n} \leq \frac{1}{\min(a_{k_j},b_{k_j})} \leq\frac{1}{f_j}$$
for any $l > 0$, and therefore
$$\lim_{n\rightarrow \infty} \frac{\max \{j: d_j \leq n\} }{n} = 0.$$
This completes the proof of the lemma.
\endproof


\section{Sparse representations by multilinear forms}
\label{sparse}

In this section we are using the setup and notation of~\cite{ab-lf}. Let $1 \leq d < n$ be integers, and let $\J(n,d)$ be the set of indexing subsets as defined in~\eqref{index_set}. For each indexing set $I = \{i_1,\dots,i_d\} \in \J(n,d)$ with $1 \le i_1 < \dots < i_d \le n$, we define the monomial $x_I$ in the variables $x_{i_1},\dots,x_{i_d}$ out of the $n$ variables $x_1,\dots,x_n$ as $x_I := x_{i_1} \cdots x_{i_d}$. An {\it integer multilinear $(n,d)$-form} is a polynomial of the form
$$F(x_1,\dots,x_n) = \sum_{I \in \J(n,d)} f_I x_I \in \zed[x_1,\dots,x_n].$$
Such an $F$ is a homogeneous polynomial in $n$ variables of degree $d$ which has degree $1$ in each of the variables $x_1,\dots,x_n$. From here on we also assume that the form $F$ is {\it coprime}, meaning that $\gcd(f_I : I \in \J(n,d)) = 1$. We say that an integer $b$ has a representation by $F$ if there exists a nonzero vector $\ba \in \zed^n$ such that $F(\ba) = b$. In this case, we also say that $F$ represents $b$ by $\ba$.

We call a nonzero integer vector $\ba$ {\it $k$-sparse} for some $1 \leq k \leq n$ if $\ba$ has no more than $k$ nonzero coordinates. We say that an integer $b$ has a $k$-sparse representation by $F$ (or $F$ represents $b$ $k$-sparsely) if there exists a $k$-sparse nonzero vector $\ba \in \zed^n$ such that $F(\ba) = b$. In~\cite{ab-lf} some results on existence of small-norm representations of integers by multilinear forms have been established. Here we build on the results of~\cite{ab-lf} to make some simple observations on the existence of integer sparse representations of an arbitrary integer $b$ by~$F$.

\begin{lem} \label{coeff-1} Suppose that one of the coefficients of $F$ is equal to $\pm 1$. Then $F$ represents $b$ by a $d$-sparse vector $\ba \in \zed^n$ with one coordinate equal to $\pm b$ and the rest of the nonzero coordinates equal to~$1$.
\end{lem}

\proof
Let $\pm x_{i_1} \cdots x_{i_d}$ be the monomial with coefficient $\pm 1$ in $F$. We can assume that $i_1 = 1,\dots, i_d = d$ without loss of generality. Since $F$ is a form of degree $d$, every other monomial in $F$ must contain a variable $x_j$ with $j \neq 1,\dots,d$. Then setting $x_j = 0$ for every $j \neq 1,\dots,d$ in $F$, we obtain 
$$F(x_1,\dots,x_d,0,\dots,0) = \pm x_1 \cdots x_d.$$
Setting $x_1 = \pm b$ as needed, and $x_2 = \dots = x_d = 1$, we then obtain 
$$F(\pm b,1,\dots,1,0,\dots,0) = b.$$
\endproof

\begin{lem} \label{ml_sparse-1} Suppose that $F$ has all nonzero pairwise coprime coefficients, none of which are equal to~$\pm 1$. Then $F$ represents every integer $k$-sparsely if and only if $k \geq d+1$. If this is the case, then for every $b \in \zed$ there exists a $k$-sparse vector $\bz \in \zed^n$ such that $F(\bz) = 0$ and
\begin{equation}
\label{z_bnd}
|\bz| \leq |b| \left( 2 |F| \right)^{\nu_d},
\end{equation}
where the exponent $\nu_{d} = \sum_{k=0}^{d} \frac{d!}{k!}$ as in equation (1) of~\cite{ab-lf}, and $|F|$ stands for the maximum of absolute values of the coefficients of~$F$.
\end{lem}

\proof
By Theorem~1.1 of~\cite{ab-lf} we know that $b$ has a representation by $F$, hence we need to address specifically the existence of sparse representations. We assume that the coefficient $f_I \neq 0,\pm 1$ for every $I \in \J(n,d)$, and the total number of coefficients is $|\J(n,d)| = \binom{n}{d}$. Notice that each variable $x_j$ is contained in $\binom{n-1}{d-1}$ monomials, hence setting $x_1 = 0$ annihilates $\binom{n-1}{d-1}$ monomials. Then setting $x_2 = 0$ annihilates $\binom{n-2}{d-1}$ of the remaining monomials. Continuing in the same manner, we see that setting $n-k$ variables equal to $0$, we annihilate
\begin{equation}
\label{sum1}
\sum_{j=1}^{n-k} \binom{n-j}{d-1} = \sum_{j=k}^{n-1} \binom{j}{d-1} = \sum_{j=0}^{n-1} \binom{j}{d-1} - \sum_{j=0}^{k-1} \binom{j}{d-1}
\end{equation}
monomials. A standard combinatorial formula asserts that
$$\sum_{j=0}^n \binom{j}{m} = \binom{n+1}{m+1}.$$
Applying this formula to~\eqref{sum1}, we obtain
$$\binom{n}{d} - \binom{k}{d}$$
annihilated monomials when $n-k$ variables are set to $0$. Since the total number of monomials is $\binom{n}{d}$, this means that the number of remaining monomials is 
$$\binom{n}{d} - \left( \binom{n}{d} - \binom{k}{d} \right) = \binom{k}{d}.$$
Notice that this number is $\geq 2$ if and only if $k \geq d+1$. Let $J_k = \{ i_1,\dots,i_k \}$ be some collection of $k$ indices in $[n]$ and let $F_{J_k}$ be the section of $F$ obtained by setting the $n-k$ variables in $[n] \setminus J_k$ equal to~$0$. Then $F_{J_k}$ is a coprime integer multilinear $(k,d)$-form with nonzero pairwise coprime coefficients, and hence represents all the integers by Theorem~1.1 of~\cite{ab-lf}.

On the other hand, if $k < d+1$ then
$$\binom{k}{d} \leq \binom{d}{d} = 1,$$
which means that setting $n-k$ coordinates equal to $0$ leaves at most one nonzero monomial, which has coefficient $\neq \pm 1$, and thus cannot represent all the integers: every integer it represents has to be divisible by its coefficient. Therefore $F$ cannot $k$-sparsely represent all the integers if $k < d+1$.

Finally, assume $k \geq d+1$ and let $b \in \zed$. In this case Theorem~1.1 of~\cite{ab-lf} guarantees the existence of a vector $\ba \in \zed^k$ such that $F_{J_k}(\ba) = b$ and 
$$|\ba| \leq |b| \left( 2 |F_{J_k}| \right)^{\nu_d},$$
since $|F_{J_k}| \leq |F|$. Take $\bz \in \zed^n$ to be the vector with coordinates
$$z_{i_j} = \left\{ \begin{array}{ll}
a_j & \mbox{if $i_j \in J_k$} \\
0 & \mbox{if $i_j \notin J_k$},
\end{array}
\right.$$
i.e. $\bz$ is formed from $\ba$ by setting every coordinate whose index is not in $J_k$ to $0$. Then $|\bz| = |\ba|$ and $F(\bz) = F_{J_l}(\ba) = b$, and so we have~\eqref{z_bnd}.
\endproof


\bibliographystyle{plain}  

\end{document}